\documentclass[12pt]{article}
\usepackage{amssymb}
\usepackage{amsfonts}
\usepackage{graphicx}
\usepackage{amsmath}
\usepackage{booktabs}
\usepackage{multirow}

\newtheorem{algorithm}{Algorithm}

\newtheorem{remark}{Remark}

\begin{document}
\date{}
\title{A Gauss Laguerre approach for the resolvent of fractional powers}
\author{Eleonora Denich\thanks{%
Dipartimento di Matematica e Geoscienze, Universit\`{a} di Trieste, Trieste, Italy, 
eleonora.denich@phd.units.it} \and Laura Grazia Dolce \thanks{Dipartimento di Matematica e Geoscienze, Universit\`{a} di Trieste, Trieste, Italy, 
lauragrazia.dolce@studenti.units.it} \and Paolo Novati \thanks{%
Dipartimento di Matematica e Geoscienze, Universit\`{a} di Trieste, Trieste, Italy, 
novati@units.it} }
\maketitle

\begin{abstract}
This paper introduces a very fast method for the computation of the resolvent of fractional powers of operators. The analysis is kept in the continuous setting of (potentially unbounded) self adjoint positive operators in Hilbert spaces. The method is based on the Gauss-Laguerre rule, exploiting a particular integral representation of the resolvent. We provide sharp error estimates that can be used to a priori select the number of nodes to achieve a prescribed tolerance.
\end{abstract}


\section{Introduction}

Let $\mathcal{H}$ be a Hilbert space with inner product $\langle \cdot,\cdot \rangle $ and induced norm $\Vert \cdot \Vert $. In this work we are interested in the computation of the resolvent of the fractional power
\begin{equation*}
\mathcal{R}_{h,\alpha}(\mathcal{L} )=(I+h\mathcal{L}^{\alpha })^{-1},
\end{equation*}%
where $h>0$, $0<\alpha <1$, $I$ is the identity operator and $\mathcal{L}:\mathcal{H}\rightarrow \mathcal{H%
}$ is a self-adjoint, positive operator. The problem finds application in the solution of fractional in space parabolic equations like 
\begin{equation*}
\frac{\partial U}{\partial t}=\Delta _{\alpha }U+F,  
\end{equation*}
where $\Delta _{\alpha }$ denotes the fractional Laplacian that can be defined following the rule of the operational calculus as $\Delta _{\alpha}=-(-\Delta )^{\alpha }$, in which $\Delta $ is the Laplace operator
with suitable boundary conditions. In the above equation $F$ represents a generic forcing term. Starting from the spectral decomposition of $\Delta $, the operation $(-\Delta )^{\alpha }$ is obtained by simply rising to the power $\alpha $ the eigenvalues of $-\Delta $. Denoting by $A\in \mathbb{R}^{N\times N}$ a general symmetric positive definite discretization of $-\Delta $, the use of an implicit time stepping procedure leads to the computation of one or more matrix functions of the type $\mathcal{R}_{h,\alpha}(A)$, where $h>0$ is a parameter that typically depends on the time step and on the coefficients of the integrator.

Going back to a generic $\mathcal{L} \colon \mathcal{H} \rightarrow \mathcal{H}$, clearly the numerical approximation of $\mathcal{R}_{h,\alpha}(\mathcal{L} )$ requires the discretization of $\mathcal{L}$. Anyway, since we want to keep the analysis independent of the type and the sharpness of the
discretization, we prefer to work in the infinite dimensional setting of the operator $\mathcal{L}$, with spectrum $\sigma (\mathcal{L})$ contained in $[a,+\infty )$, $a>0$, so potentially unbounded. Since we do not
introduce any restriction on the magnitude of $h$, without loss of generality throughout the paper we assume for simplicity $a=1$. 

The numerical approach considered in this work employs a particular integral representation of the function $\mathcal{R}_{h,\alpha}(\lambda )=(1+h\lambda^{\alpha })^{-1}$. With some manipulations, the representation allows to write $\mathcal{R}_{h,\alpha}(\lambda )$ as the sum of two
integrals, that is,
\begin{equation} \label{I}
\mathcal{R}_{h,\alpha}(\lambda )=\frac{\sin (\alpha \pi )}{\alpha
\pi }\left( I^{(1)}(\lambda )+I^{(2)}(\lambda
)\right),
\end{equation}
with 
\begin{eqnarray}
I^{(1)}(\lambda ) &=&\int_{0}^{+\infty }\frac{e^{-x}}{\left( 1+e^{-\frac{x}{%
\alpha }}h^{\frac{1}{\alpha }}\lambda \right) \left( e^{-2x}+2e^{-x}\cos
(\alpha \pi )+1\right) }dx, \label{I1}\\
I^{(2)}(\lambda ) &=&\int_{0}^{+\infty } \frac{\alpha (\alpha+1)^{-1} e^{-x}}{ \left( e^{-\frac{x}{%
\alpha +1}}+h^{\frac{1}{\alpha }}\lambda \right) \left( e^{-\frac{2\alpha x}{%
\alpha +1}}+2e^{-\frac{\alpha x}{\alpha +1}}\cos (\alpha \pi )+1\right) }dx, \label{I2}
\end{eqnarray}
that we evaluate by using the $n$-point Gauss-Laguerre rule. This formula leads to a rational approximation $R_{2n-1,2n}(\lambda)\cong \mathcal{R}_{h,\alpha}(\lambda )$, where $R_{2n-1,2n}=p_{2n-1}/q_{2n}$, with $p_{2n-1}\in \Pi _{2n-1}$, $q_{2n}\in \Pi _{2n}$.
Here and below we use the symbol $\cong$ to indicate a generic approximation.
Finally, we thus obtain 
\begin{equation*}
\mathcal{R}_{h,\alpha}(\mathcal{L} )\cong R_{2n-1,2n}(\mathcal{L}),
\end{equation*}%
in which the degree $2n$ of the denominator also represents the number of inversions of suitable shifts of the operator $\mathcal{L}$.

Thanks to the existing error estimates based on the theory of analytic functions, we are able to derive sharp error estimates for the operator
case. This is possible because ($\sigma(\mathcal{L}) \subseteq [1, +\infty)$)
\begin{equation}
\Vert \mathcal{R}_{h,\alpha}(\mathcal{L} )-R_{k-1,k}(\mathcal{L})\Vert _{\mathcal{H}\rightarrow \mathcal{H}}\leq \max_{\lambda\in \lbrack 1,+\infty
)}|\mathcal{R}_{h,\alpha}(\lambda )-R_{2n-1,2n}(\lambda)|,  \label{2}
\end{equation}
where $\Vert \cdot \Vert _{\mathcal{H}\rightarrow \mathcal{H}}$ is the induced operator norm. The bound (\ref{2}) allows to study the error by working scalarly.
This analysis shows that, by using the same number $n$ of Laguerre points for both integrals, the error decays like
\begin{equation*}
\exp \left( - \emph{const} \cdot n^{1/3}\right).
\end{equation*}%
Nevertheless, exploiting the sharpness of the estimates we are able to remarkably improve the efficiency of the method by balancing the error contributions of the two integrals, together with a suitable truncation of the
Laguerre rule. In this way we show that the error decays like
\begin{equation*}
\exp \left( -\emph{const} \cdot q^{1/2}\right),
\end{equation*}
where $q \ll 2n$ represents the number of inversions.

We remark that the analysis given in the paper is rather complicate and we have been forced to use a number of approximations, sometimes consequence of experimental evidences. Nevertheless, the final estimates appear to be very accurate and this, somehow, justifies our choices.

The computation of the resolvent $\mathcal{R}_{h,\alpha}(\mathcal{L} )$, where $\mathcal{L}$ is an unbounded operator, has been already studied in \cite{AN20}. The particular case of $\mathcal{L}$ representing a matrix has been considered in \cite{MN}), where the shift and invert Krylov method is employed, and in \cite{ABDN}, where the rational Krylov method based on the use of zeros of the Jacobi polynomials is considered. We remark that the approach here proposed can also be employed to define the poles of a rational Krylov method. 
Since $\mathcal{R}_{h,\alpha}(\lambda) \sim \frac{1}{h} \lambda^{-\alpha}$ for $\lambda \rightarrow + \infty$ (the symbol $\sim$ denotes the asymptotic equivalence), in principle one may use the poles of a rational approximation of $\lambda^{-\alpha}$ (see e.g., \cite{AN19, BP, BLP, HLMMV, HLMMP, HLM, HLMMPbook, HMbook, V2015, V2018, V2020}) also for $\mathcal{R}_{h,\alpha}(\lambda)$. Anyway, clearly, working directly with $\mathcal{R}_{h,\alpha}(\lambda)$ allows to obtain better results.

The work is organized as follows. In Section \ref{Section2} we derive the integral representation used and describe the Gauss Laguerre method. In Section \ref{Section3} we study the error in the scalar case. In Section \ref{Section4} we extend the analysis to the operator case. Finally, in Section \ref{Section5} and \ref{Section6} we present some improvements based on the properties of the integrand functions and on the asymptotic behavior of the Laguerre weights.

\section{The Gauss-Laguerre method} \label{Section2}

Let $\lambda \in \mathbb{C}\setminus (-\infty ,0]$ and consider the Cauchy integral representation 
\begin{equation*}
\mathcal{R}_{h,\alpha}(\lambda )=\frac{1}{2\pi i}\int_{\Gamma }(z-\lambda
)^{-1}(1+hz^{\alpha })^{-1}dz,
\end{equation*}%
where $\Gamma $ is a contour in $\mathbb{C}\setminus (-\infty ,0)$ containing $\lambda $ in its interior. The first step for the construction of our method is to select $\Gamma $ as the boundary of the sector of the
complex plane with semiangle $\alpha \pi $ and vertex at the origin, that is,%
\begin{equation*}
\Gamma =\Gamma _{\alpha }:=\Gamma _{\alpha }^{+}\cup \Gamma _{\alpha }^{-},
\end{equation*}%
where $\Gamma _{\alpha }^{\pm }=\{w|w=\rho e^{\pm i\alpha \pi },\rho \geq
0\} $.
Then, the change of variable $z=\left( \frac{w}{h}\right) ^{1/\alpha}$ leads to 
\begin{equation*}
\mathcal{R}_{h,\alpha}(\lambda )=\frac{1}{2\pi i\alpha h^{\frac{1}{%
\alpha }}}\int_{\Gamma _{\alpha }}\frac{1}{\left( w^{\frac{1}{\alpha }}h^{-%
\frac{1}{\alpha }}-\lambda \right) \left( 1+w\right) }w^{\frac{1}{\alpha }%
-1}dw.
\end{equation*}%
Now, defining $w=\rho e^{i\pi \alpha }$ for $w\in \Gamma _{\alpha }^{+}$, $%
w=\rho e^{-i\pi \alpha }$ for $w\in \Gamma _{\alpha }^{-}$, and running $%
\Gamma _{\alpha }$ in the counterclockwise direction, we obtain%
\begin{align}
\mathcal{R}_{h,\alpha}(\lambda ) & =\frac{1}{2\pi i\alpha h^{\frac{1%
}{\alpha }}}\left[ -\int_{0}^{+\infty }\frac{\rho ^{\frac{1}{\alpha }e^{i\pi
}e^{i\alpha \pi }}}{\rho e^{i\pi \alpha }\left( \rho ^{\frac{1}{\alpha }%
}e^{i\pi }h^{-\frac{1}{\alpha }}-\lambda \right) \left( 1+\rho e^{+i\alpha
\pi }\right) }d\rho \right. \notag \\
& \left. +\int_{0}^{+\infty }\frac{\rho ^{\frac{1}{\alpha }e^{-i\pi
}e^{-i\alpha \pi }}}{\rho e^{-i\pi \alpha }\left( \rho ^{\frac{1}{\alpha }%
}e^{-i\pi }h^{-\frac{1}{\alpha }}-\lambda \right) \left( 1+\rho e^{-i\alpha
\pi }\right) }d\rho \right] \notag \\ 
&=\frac{\sin (\alpha \pi )}{\alpha
\pi }\int_{0}^{+\infty }\frac{\rho ^{\frac{1}{\alpha }}}{\left( \rho ^{\frac{%
1}{\alpha }}+h^{\frac{1}{\alpha }}\lambda \right) \left( 1+2\rho \cos
(\alpha \pi )+\rho ^{2}\right) }d\rho .  \label{K}
\end{align}
We remark that the integral representation (\ref{K}), with $\lambda $ replaced by a regularly accretive operator $\mathcal{L}$, was derived by Kato in \cite{Kato}.
Now, by the change of variable $\rho =e^{y}$, we obtain 
\begin{align*}
\mathcal{R}_{h,\alpha}(\lambda ) &=\frac{\sin (\alpha \pi )}{\alpha
\pi }\int_{-\infty }^{+\infty }\frac{e^{\frac{\alpha +1}{\alpha }y}}{\left(
e^{\frac{y}{\alpha }}+h^{\frac{1}{\alpha }}\lambda \right) \left( 1+2\cos
(\alpha \pi )e^{y}+e^{2y}\right) }dy \\
& =\frac{\sin (\alpha \pi )}{\alpha
\pi }\left[ \int_{0}^{+\infty }\frac{e^{-y}}{\left( 1+e^{\frac{-y}{\alpha }%
}h^{\frac{1}{\alpha }}\lambda \right) \left( e^{-2y}+2\cos (\alpha \pi
)e^{-y}+1\right) }dy\right. \\
& \left. +\int_{-\infty }^{0}\frac{e^{\frac{\alpha +1}{\alpha }y}}{\left( e^{%
\frac{y}{\alpha }}+h^{\frac{1}{\alpha }}\lambda \right) \left( 1+2\cos
(\alpha \pi )e^{y}+e^{2y}\right) }dy\right].
\end{align*}
Using for the second integral the further change of variable $x=\frac{\alpha +1}{\alpha }y,$ we finally obtain (\ref{I}), in which the integrals (see (\ref{I1}) and (\ref{I2}))   can be written as
\begin{equation*}
I^{(i)}(\lambda)=\int_{0}^{+\infty }e^{-x}f_{i}(x)dx,\quad i=1,2, 
\end{equation*}%
with 
\begin{eqnarray}
f_{1}(x) &=&\left( 1+e^{\frac{-x}{\alpha }}h^{\frac{1}{\alpha }}\lambda
\right) ^{-1}\left( e^{-2x}+2e^{-x}\cos (\alpha \pi )+1\right) ^{-1}, \label{f1}\\
f_{2}(x) &=& \frac{\alpha}{\alpha+1} \left( e^{\frac{-x}{\alpha +1}}+h^{\frac{1}{\alpha }}\lambda
\right) ^{-1}\left( 1+2\cos (\alpha \pi )e^{\frac{-\alpha x}{\alpha +1}}+e^{%
\frac{-2\alpha x}{\alpha +1}}\right) ^{-1}. \label{f2}
\end{eqnarray}
\begin{remark}
The functions $f_1$ and $f_2$ depend on $\lambda$ in the scalar case or on $\mathcal{L}$ in the operator case. Anyway, to keep the notations as simple as possible, we omit this dependency since it is always clear from the context.
\end{remark}

Remembering that we assume $\sigma (\mathcal{L})\subseteq \lbrack 1,+\infty )$, we restrict our considerations to the case $\lambda \geq 1$. In this situation, we observe that for $x\in \lbrack 0,+\infty )$ we have $0\leq
f_{i}(x)\leq K_{i}$, $i=1,2$, where 
\begin{equation} \label{K_i}
K_{1}=1 \quad {\rm and} \quad K_{2}=\frac{\alpha}{\alpha+1}h^{-\frac{1}{\alpha }}.
\end{equation}
By using the Gauss-Laguerre rule, we then approximate the resolvent (cf. (\ref{I})) by 
\begin{equation}
\mathcal{R}_{h,\alpha}(\lambda ) \cong \frac{\sin (\alpha \pi )}{\alpha \pi }%
\left[ I_{n}^{(1)}(\lambda )+I_{n}^{(2)}(\lambda )%
\right],   \label{lag}
\end{equation}%
where 
\begin{equation*}
I_{n}^{(i)}(\lambda ):=\sum_{j=1}^{n}w_{j}^{(n)}f_{i}(x_{j}^{(n)}),
\end{equation*}
in which $x_{j}^{(n)}$ and $w_{j}^{(n)}$, for $j=1,...,n$, are respectively the nodes and the weights of the $n$-point Laguerre rule. 
In the sequel we denote by  
\begin{equation} \label{e}
e_n^{(i)} =I^{(i)}(\lambda)-I_n^{(i)}(\lambda), \quad i=1,2,
\end{equation}
the corresponding errors.
Clearly, (\ref{lag}) is
actually a rational approximation of type 
\begin{equation*}
\mathcal{R}_{h,\alpha}(\lambda ) \cong R_{2n-1,2n}(\lambda) = \frac{p_{2n-1}(\lambda )}{q_{2n}(\lambda )},\quad p_{2n-1}\in \Pi _{2n-1}, \; q_{2n}\in \Pi _{2n}.
\end{equation*}

\section{Error estimates in the scalar case} \label{Section3}

Let $f$ be a generic function, analytic in a region of the complex plane containing the positive real axis. Let 
\begin{equation*}
I(f)=\int_{0}^{+\infty }e^{-x}f(x)dx
\end{equation*}%
and $I_{n}(f)$ be the $n$-point Gauss-Laguerre rule. 
In order to estimate $e_{n}(f)=I(f)-I_{n}(f)$, we employ the theory introduced in \cite{Barrett}.
Let us define for any $\mathcal{S}>1$ the parabola $\Gamma _{\mathcal{S}}$ in the complex plane given by 
\begin{equation}
\mathfrak{Re}(\sqrt{-z})=\ln (\mathcal{S}).  \label{par}
\end{equation}%
$\Gamma _{\mathcal{S}}$ is symmetric with respect to the real axis, it has convexity oriented towards the positive real axis and vertex in $-(\ln (\mathcal{S}))^{2}$. By writing $z=a+ib$, equation (\ref{par}) becomes 
\begin{equation*}
\sqrt{\frac{\sqrt{a^{2}+b^{2}}-a}{2}}=\ln (\mathcal{S}),
\end{equation*}
and then we can rewrite the parabola in Cartesian coordinates as 
\begin{equation}
a=\left( b^{2}-4(\ln (\mathcal{S}))^{4}\right) \frac{1}{4(\ln (\mathcal{S}))^{2}}.  \label{cart}
\end{equation}
Observe that for $\mathcal{S}\rightarrow 1$ the parabola degenerates to the positive real axis.

Suppose that for a certain $\mathcal{S}$ the function $f$ is analytic on and within $\Gamma _{\mathcal{S}}$ except for a pair of simple poles $z_{0}$ and its conjugate $\overline{z}_{0}$. Then, the error can be estimated as follows (see \cite{Barrett}) 
\begin{equation}
e_{n}(f)\cong -4\pi \mathfrak{Re}\left\{ r(f;z_{0})e^{-z_{0}}\left( e^{\sqrt{%
-z_{0}}}\right) ^{-2\sqrt{\overline{n}}}\right\} , \label{berr}
\end{equation}
where $\overline{n}=4n+2$ and $r(f;z_{0})$ is the residue of $f$ at $z_{0}$.
Let $\mathcal{S}_{0}>1$ be such that $z_{0}$ and $\overline{z}_{0}$ belong
to $\Gamma _{\mathcal{S}_{0}}$. Then, by (\ref{par}) and (\ref{berr}) we have that
\begin{equation}
\left\vert e_{n}(f)\right\vert \cong const \cdot\mathcal{S}_{0}^{-2\sqrt{\overline{n}}}.
\label{s0}
\end{equation}

\subsection{Poles and residues}

In order to employ the above theory to estimate the error of approximation (\ref{lag}), we have to study poles and residues of the functions $f_{1}(z)$ and $f_{2}(z)$ (see (\ref{f1}) and (\ref{f2})). For both functions, the key point is to understand which are the poles closest to $[0,+\infty )$, where the distance is expressed by $\mathcal{S}_{0}$ (see (\ref{s0})). For simplicity, we define the auxiliary
functions
\begin{align*}
a^{(I)}(z)& =1+e^{\frac{-z}{\alpha }}h^{\frac{1}{\alpha }}\lambda , \\
a^{(II)}(z)& =e^{-2z}+2e^{-z}\cos (\alpha \pi )+1, \\
a^{(III)}(z)& =e^{\frac{-z}{\alpha +1}}+h^{\frac{1}{\alpha }}\lambda , \\
a^{(IV)}(z)& =1+2\cos (\alpha \pi )e^{\frac{-\alpha z}{\alpha +1}}+e^{\frac{%
-2\alpha z}{\alpha +1}},
\end{align*}
so that%
\begin{equation*}
f_{1}(z)=\left[ a^{(I)}(z)a^{(II)}(z)\right] ^{-1},\quad f_{2}(z)=\frac{\alpha}{\alpha+1} \left[
a^{(III)}(z)a^{(IV)}(z)\right] ^{-1}.
\end{equation*}

\subsubsection{First integral} \label{section3.1.1}

The poles of $f_{1}$ are obtained by solving $%
a^{(I)}(z)a^{(II)}(z)=0$. Starting from $a^{(I)}(z)=0$ we obtain the set
\begin{equation*}
z_{k}^{(I)}=\alpha \ln (h^{\frac{1}{\alpha }}\lambda )+i(2k+1)\alpha \pi
,\quad k\in \mathbb{Z}.
\end{equation*}%
The closest to $[0,+\infty )$ are therefore%
\begin{equation}
z_{0}^{(I)}=\alpha \ln (h^{\frac{1}{\alpha }}\lambda )+i\alpha \pi 
\label{p1}
\end{equation}%
and its conjugate, obtained for $k=0$ and $k=-1$. As for the poles arising
from $a^{(II)}(z)$ we obtain 
\begin{equation*}
z_{k}^{(II)}=i(2k+1\pm \alpha )\pi ,\quad k\in \mathbb{Z},
\end{equation*}%
and now the closest to the real axis are 
\begin{equation}
z_{0}^{(II)}=i(1-\alpha )\pi   \label{p2}
\end{equation}%
and its conjugate. In order to estimate the error $e_{n}^{(1)}(\lambda )$ (see (\ref{e})), let $\Gamma _{\mathcal{S}_{0}^{(II)}}$ be the parabola passing through the pole $z_{0}^{(II)}$. If it contains $z_{0}^{(I)}$ in its interior we use the error formula (\ref{berr}) with $z_0=z_{0}^{(I)}$, otherwise with $z_0=z_{0}^{(II)}$. By using the Cartesian expression (\ref{cart}), it is rather easy to
demonstrate that $z_{0}^{(I)}$ is inside $\Gamma _{\mathcal{S}_{0}^{(II)}}$ for 
\begin{equation}
\ln \left(h^{\frac{1}{\alpha }}\lambda \right)>\frac{(2\alpha -1)\pi }{2\alpha
(1-\alpha )}.  \label{c1}
\end{equation}
Defining
\begin{equation}
\overline{\lambda }=\exp \left( \frac{(2\alpha -1)\pi }{2\alpha (1-\alpha )}%
\right) h^{-\frac{1}{\alpha }},  \label{lb}
\end{equation}
we simply write $\lambda >\overline{\lambda }$ to express condition (\ref{c1}). We
remark that, depending on $h$ and $\alpha $, this condition may be verified for each $\lambda \geq 1$, or eventually, only for $\lambda $ sufficiently large.

As for the residues of $f_{1}$ at $z_{0}^{(I)}$ and $z_{0}^{(II)}$, after some computations one finds
\begin{eqnarray*}
r(f_{1};z_{0}^{(I)}) &=&\frac{1}{a^{(II)}(z_{0}^{(I)})}\lim_{z\rightarrow
z_{0}^{(I)}}\frac{z-z_{0}^{(I)}}{a^{(I)}(z)}  \notag \\
&=&\frac{\alpha h^{2}\lambda ^{2\alpha }}{e^{-2i\alpha \pi }+2h\lambda
^{\alpha }\cos (\alpha \pi )e^{-i\alpha \pi }+h^{2}\lambda ^{2\alpha }},
\end{eqnarray*}%
and%
\begin{eqnarray*}
r(f_{1};z_{0}^{(II)}) &=&\frac{1}{a^{(I)}(z_{0}^{(II)})}\lim_{z\rightarrow
z_{0}^{(II)}}\frac{z-z_{0}^{(II)}}{a^{(II)}(z)}  \notag \\
&=&\frac{ie^{-i\alpha \pi }}{2\sin (\alpha \pi )\left( 1-e^{-i\frac{\pi }{%
\alpha }}h^{\frac{1}{\alpha }}\lambda \right) }. 
\end{eqnarray*}

\subsubsection{Second integral}

Following the same steps of Section \ref{section3.1.1}, we obtain the poles of $f_{2}$ by solving $%
a^{(III)}(z)a^{(IV)}(z)=0$. Starting from $a^{(III)}(z)=0$ we find the set%
\begin{equation*}
z_{k}^{(III)}=-(\alpha +1)\ln (h^{\frac{1}{\alpha }}\lambda )+i(2k+1)(\alpha
+1)\pi ,\quad k\in \mathbb{Z},
\end{equation*}%
and the closest to the real axis are 
\begin{equation}
z_{0}^{(III)}=-(\alpha +1)\ln (h^{\frac{1}{\alpha }}\lambda )+i(\alpha
+1)\pi   \label{p3}
\end{equation}%
and its conjugate. By solving $a^{(IV)}(z)=0$ we obtain
\begin{equation*}
z_{k}^{(IV)}=i(2k+1\pm \alpha )\frac{\alpha +1}{\alpha }\pi ,\quad k\in 
\mathbb{Z},
\end{equation*}%
and now the closest are
\begin{equation}
z_{0}^{(IV)}=i\frac{(1-\alpha )(\alpha +1)}{\alpha }\pi  , \label{p4}
\end{equation}
and its conjugate.
Analogously to the first integral, we have to compare the parabolas passing through $z_{0}^{(III)}$ and $z_{0}^{(IV)}$. With some computations it is not
difficult to see that $z_{0}^{(III)}$ is inside the parabola $\Gamma _{\mathcal{S}_{0}^{(IV)}}$ passing though $z_{0}^{(IV)}$ if
\begin{equation*}
\ln \left(h^{\frac{1}{\alpha }}\lambda \right)<-\frac{(2\alpha -1)\pi }{2\alpha
(1-\alpha )}. 
\end{equation*}
It may happen that the above condition is not verified for any $\lambda \geq 1$. Hence, it must be replaced by 
\begin{equation}
1\leq \lambda <  \bar{\bar{\lambda}}, \quad {\rm where} \quad \bar{\bar{\lambda}}=\max \left\lbrace 1, e^{-\frac{(2\alpha -1)\pi }{2\alpha (1-\alpha )}}h^{-\frac{1}{\alpha }} \right\rbrace.  \label{c3}
\end{equation}%
If $\bar{\bar{\lambda}}=1$, then the pole to consider is
always $z_{0}^{(IV)}$.

For what concerns the corresponding residues, we obtain%
\begin{eqnarray*}
r(f_{2};z_{0}^{(III)}) &=& \frac{\alpha}{(\alpha +1) a^{(IV)}(z_{0}^{(III)})}%
\lim_{z\rightarrow z_{0}^{(III)}}\frac{z-z_{0}^{(III)}}{a^{(III)}(z)}  \notag
\\
&=&\frac{\alpha}{h^{\frac{1}{\alpha }}\lambda \left( 1+2\cos (\alpha \pi
)(-1)^{\alpha }h\lambda ^{\alpha }+(-1)^{2\alpha }h^{2}\lambda ^{2\alpha
}\right) } ,
\end{eqnarray*}%
and%
\begin{eqnarray*}
r(f_{2};z_{0}^{(IV)}) &=&\frac{\alpha}{(\alpha+1)a^{(III)}(z_{0}^{(IV)})}\lim_{z\rightarrow
z_{0}^{(IV)}}\frac{z-z_{0}^{(IV)}}{a^{(IV)}(z)}  \notag \\
&=&-\frac{ie^{i\alpha \pi }}{2 \sin (\alpha \pi )\left( e^{\frac{i(1-\alpha )\pi }{\alpha }}+h^{\frac{1}{\alpha }%
}\lambda \right) }.  
\end{eqnarray*}

\subsection{The final estimates} \label{Section3.2}

\begin{figure}
\begin{center}
\includegraphics[scale=0.35]{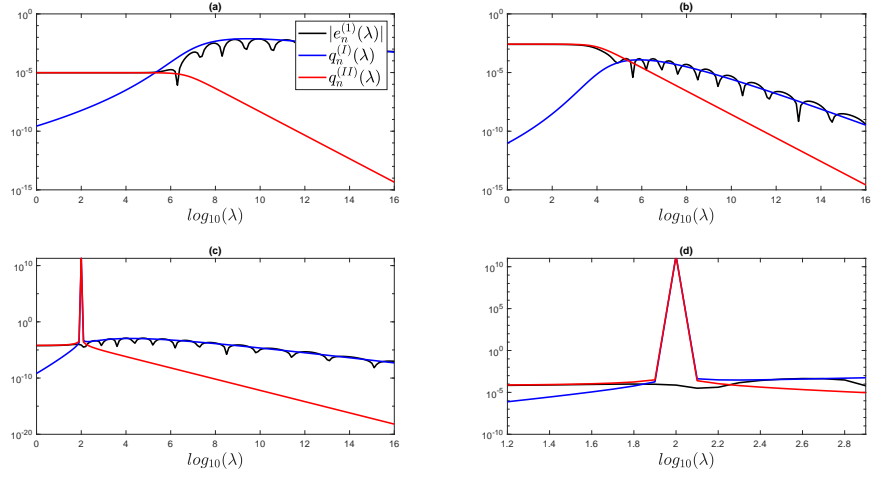}
\end{center}
\caption{Behavior of the functions $\left\vert e_n^{(1)}(\lambda)\right\vert$, $q_n^{(I)}(\lambda)$, $q_n^{(II)}(\lambda)$ for $\alpha=0.3$, $h=1e-2$ (a), $\alpha=0.75$, $h=1e-3$ (b), $\alpha=0.5$, $h=1e-1$ (c). Picture (d) is a zoom of (c) around $\bar{\lambda}$. In all cases $n=30$.}  \label{Figura1}
\end{figure}

In order to apply (\ref{berr}) it remains to evaluate the terms $\exp\left(-z_{0}^{(\cdot)}\right)$ and $\exp \left(\sqrt{-z_{0}^{(\cdot)}}\right)$, where $z_{0}^{(\cdot)}$ represents one of the four poles to be considered. By (\ref{p1}), (\ref{p2}), (\ref{p3}), (\ref{p4}), we immediately have
\begin{eqnarray*}
\exp (-z_{0}^{(I)}) &=&\frac{e^{-i\alpha \pi }}{h\lambda ^{\alpha }}, \\
\exp (-z_{0}^{(II)}) &=&-e^{i\alpha \pi }, \\
\exp (-z_{0}^{(III)}) &=&(-1)^{\alpha +1}h^{\frac{\alpha +1}{\alpha }%
}\lambda ^{\alpha +1}, \\
\exp (-z_{0}^{(IV)}) &=&e^{\frac{i(1-\alpha )(\alpha +1)\pi }{\alpha }}.
\end{eqnarray*}%
As for the quantity $\sqrt{-z_{0}^{(\cdot)}}$, using the relation
\begin{equation*}
\sqrt{-iw}=\frac{\sqrt{2}}{2}\left( 1-i\right) \sqrt{w},\quad w\geq 0,
\end{equation*}%
we obtain 
\begin{align*}
\sqrt{-z_{0}^{(I)}}& =\sqrt{\frac{\alpha }{2}}\left( \gamma ^{-}(\lambda
)-i\gamma ^{+}(\lambda )\right) , \\
\sqrt{-z_{0}^{(II)}}& =\frac{\sqrt{2}}{2}\left( 1-i\right) \sqrt{(1-\alpha
)\pi }, \\
\sqrt{-z_{0}^{(III)}}& =\sqrt{\frac{\alpha +1}{2}}\left( \gamma ^{+}(\lambda
)+i\gamma ^{-}(\lambda )\right) , \\
\sqrt{-z_{0}^{(IV)}}& =\frac{\sqrt{2}}{2}\left( 1-i\right) \sqrt{\frac{%
(1-\alpha )(\alpha +1)}{\alpha }\pi },
\end{align*}
where 
\begin{equation}
\gamma ^{\pm }(\lambda )=\sqrt{\sqrt{(\ln (h^{\frac{1}{\alpha }}\lambda
))^{2}+\pi ^{2}}\pm \ln (h^{\frac{1}{\alpha }}\lambda )}.  \label{gamma}
\end{equation}
Using the above results and the residues of previous section in (\ref{berr}), we are finally able to write the error estimates for both integrals. For the first one we have%
\begin{equation*}
e_{n}^{(1)}(\lambda )\cong \left\{ 
\begin{array}{c}
-4\pi \mathfrak{Re}\left( \frac{\alpha e^{-i\alpha \pi }h\lambda ^{\alpha
}e^{-\sqrt{2\alpha \overline{n}}\left( \gamma ^{-}(\lambda )+i\gamma
^{+}(\lambda )\right) }}{e^{-2i\alpha \pi }+2h\lambda ^{\alpha }\cos (\alpha
\pi )e^{-i\alpha \pi }+h^{2}\lambda ^{2\alpha }}\right) \quad \text{if }%
\lambda >\overline{\lambda } \\ 
-2\pi \mathfrak{Re}\left( \frac{-ie^{-\left( 1-i\right) \sqrt{2(1-\alpha
)\pi \overline{n}}}}{\sin (\alpha \pi )(1-e^{-i\frac{\pi }{\alpha }}h^{\frac{%
1}{\alpha }}\lambda )}\right) \quad \text{if } 1< \lambda <\overline{\lambda }
\end{array}%
\right. ,
\end{equation*}
($\bar{\lambda}$ defined in (\ref{lb})) and therefore, by taking the modulus,
\begin{equation*}
\left\vert e_{n}^{(1)}(\lambda )\right\vert \cong \left\{ 
\begin{array}{c}
q_{n}^{(I)}(\lambda )\quad \text{if }\lambda >\overline{\lambda }, \\ 
q_{n}^{(II)}(\lambda )\quad \text{if }1< \lambda <\overline{\lambda }
\end{array}%
\right. ,
\end{equation*}%
where
\begin{eqnarray}
q_{n}^{(I)}(\lambda ) &=& \frac{4\pi\alpha h\lambda ^{\alpha }e^{-\sqrt{%
2\alpha \overline{n}}\gamma ^{-}(\lambda )}}{\left\vert e^{-2i\alpha \pi
}+2h\lambda ^{\alpha }\cos (\alpha \pi )e^{-i\alpha \pi }+h^{2}\lambda
^{2\alpha }\right\vert }, \label{qI} \\
q_{n}^{(II)}(\lambda ) &=& \frac{2\pi e^{-\sqrt{2(1-\alpha )\pi 
\overline{n}}}}{\sin (\alpha \pi )\left\vert 1-e^{-i\frac{\pi }{\alpha }}h^{%
\frac{1}{\alpha }}\lambda \right\vert }. \label{qII} 
\end{eqnarray}
We refer to Figure \ref{Figura1} for some experiments, where $n=30$ is fixed and we consider, for different values of $\alpha$ and $h$, what happens for $\lambda \in \left[1, 10^{16} \right]$.
\begin{remark}
From a theoretical point of view, the estimates (\ref{qI}) and (\ref{qII}) should not work properly for $\lambda$ close to $\bar{\lambda}$ and are incorrect for $\lambda=\bar{\lambda}$. In this situation the poles $z_0^{(I)}$ and $z_0^{(II)}$ are such that the corresponding parabolas $\Gamma_{\mathcal{S}_0^{(I)}}$, $\Gamma_{\mathcal{S}_0^{(II)}}$ overlap and then the analysis should consider the contribution of both poles in the computation of the residues. 
Nevertheless, as clearly shown in Figure \ref{Figura1}a, \ref{Figura1}b, the true error is very well approximated by $q_n^{(\cdot)}(\lambda)$ in the whole interval.
In the particular case of $\lambda=\bar{\lambda}$ and $\alpha=1/2$ we have that $z_0^{(I)}=z_0^{(II)}$ and then formula (\ref{berr}) does not work, because we have a double pole (see the peak in Figure \ref{Figura1}c). 
In order to undertake this situation without weighing down the theory in Section \ref{Section4}, we employ some  simplifications that come from the evidence that, unlike the estimates, the true error does not explodes (see Figure \ref{Figura1}d, representing a zoom of Figure \ref{Figura1}c).
The same considerations hold true for the second integral whose final results are given below (see (\ref{qIII}), (\ref{qIV})).
\end{remark}

\begin{figure}
\begin{center}
\includegraphics[scale=0.35]{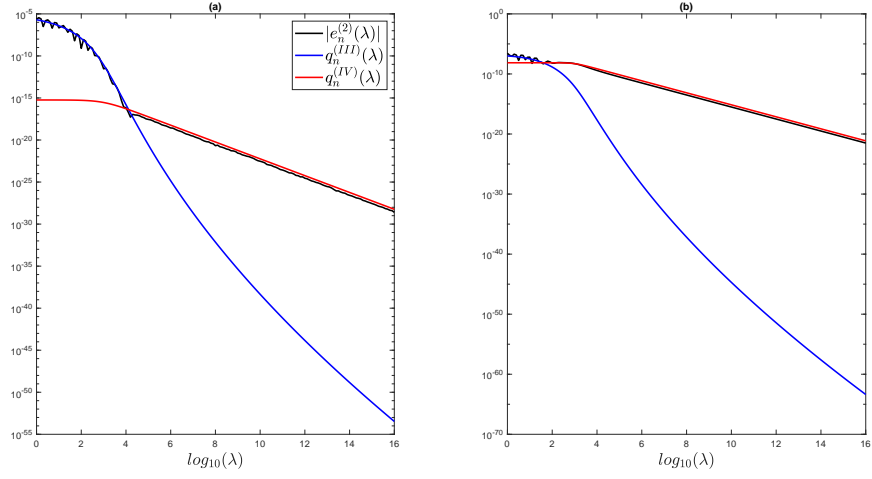}
\end{center}
\caption{Behavior of the functions $| e_n^{(2)}(\lambda)|$, $q_n^{(III)}(\lambda)$, $q_n^{(IV)}(\lambda)$ for $\alpha=1/3$, $h=1e-1$ (a) and $\alpha=2/3$, $h=1e-2$ (b). In both cases $n=20$.}  \label{Figura2}
\end{figure}

For what concerns the error of the second integral, using (\ref{berr}) and the previous results, we obtain%
\begin{equation*}
e_{n}^{(2)}(\lambda )\cong \left\{ 
\begin{array}{c}
4\pi \mathfrak{Re}\left( \frac{\alpha(-1)^{\alpha }h\lambda ^{\alpha }e^{-\sqrt{2(\alpha +1)\overline{n}}\left( \gamma
^{+}(\lambda )+i\gamma ^{-}(\lambda )\right) }}{%
1+2\cos (\alpha \pi )(-1)^{\alpha }h\lambda ^{\alpha }+(-1)^{2\alpha
}h^{2}\lambda ^{2\alpha }}\right) \quad \text{if }1\leq
\lambda <\overline{\overline{\lambda }} \\ 
2\pi \mathfrak{Re}\left\{ \frac{ie^{i\alpha \pi }e^{\frac{%
i(1-\alpha )(\alpha +1)\pi }{\alpha }}}{\sin (\alpha \pi )\left( e^{%
\frac{i(1-\alpha )\pi }{\alpha }}+h^{\frac{1}{\alpha }}\lambda \right) }%
e^{-(1-i)\left( \sqrt{\frac{2(1-\alpha )(\alpha +1)\pi }{\alpha }\overline{n}%
}\right) }\right\} \quad \text{if }\lambda > \overline{\overline{\lambda }%
}%
\end{array}%
\right. ,
\end{equation*}%
where $\overline{\overline{\lambda }}$ is defined in (\ref{c3}). Thus we have
\begin{equation*}
\left\vert e_{n}^{(2)}(\lambda )\right\vert \cong \left\{ 
\begin{array}{c}
q_{n}^{(III)}(\lambda )\quad \text{if }1\leq \lambda <\overline{\overline{%
\lambda }} \\ 
q_{n}^{(IV)}(\lambda )\quad \text{if }\lambda > \overline{\overline{%
\lambda }}%
\end{array}%
\right. ,
\end{equation*}%
where%
\begin{eqnarray}
q_{n}^{(III)}(\lambda ) &=& \frac{4\pi\alpha h\lambda ^{\alpha }e^{-%
\sqrt{2(\alpha +1)\overline{n}}\gamma ^{+}(\lambda )}}{\left\vert 1+2\cos
(\alpha \pi )(-1)^{\alpha }h\lambda ^{\alpha }+(-1)^{2\alpha }h^{2}\lambda
^{2\alpha }\right\vert }, \label{qIII} \\
q_{n}^{(IV)}(\lambda ) &=& \frac{2\pi e^{-\sqrt{\frac{2(1-\alpha
)(\alpha +1)\pi }{\alpha }\overline{n}}}}{\sin (\alpha \pi
)\left\vert e^{\frac{i(1-\alpha )\pi }{\alpha }}+h^{\frac{1}{\alpha }%
}\lambda \right\vert }. \label{qIV}
\end{eqnarray}
We refer to Figure \ref{Figura2} for a couple of experiments.

\section{Error estimate for the operator} \label{Section4}

By using (\ref{I}), (\ref{2}) and the results of the previous section, the idea now is to estimate the error of the method as follows
\begin{align*}
\Vert \mathcal{R}_{h,\alpha}(\mathcal{L} )-R_{k-1,k}(\mathcal{L})\Vert _{\mathcal{H}\rightarrow \mathcal{H}} &\cong \frac{\sin (\alpha \pi)}{\alpha \pi} \times \\\Bigg(\max \left\lbrace \max_{\lambda \geq \bar{\lambda}} q_n^{(I)}(\lambda), \max_{1 \leq \lambda \leq \bar{\lambda}} q_n^{(II)}(\lambda) \right\rbrace &+ \max \left\lbrace \max_{1 \leq \lambda \leq \bar{\bar{\lambda}}} q_n^{(III)}(\lambda), \max_{\lambda \geq \bar{\bar{\lambda}}} q_n^{(IV)}(\lambda) \right\rbrace \Bigg).
\end{align*}
The problem is then reduced to the evaluation of the maximum of the functions $q_n^{(\cdot)}$. Since these functions are not very simple to handle, we are forced to use some further approximations. Note that, in the above formula, we have included the boundaries $\bar{\lambda}$, $\bar{\bar{\lambda}}$ and then it does not work for $\alpha=1/2$ (see Figure \ref{Figura1}c, \ref{Figura1}d). Nevertheless, the approximation used will allow to solve the question.

\subsection{Approximation of maxima} \label{Section4.1}

\quad $\mathbf{ q_{n}^{(I)}(\lambda )}$ - Independently of $\alpha $ and for $n$ large enough, the function reaches a relative maximum at a certain point $\lambda_n$ (Figure \ref{Figura1}), and then goes to zero for $\lambda \rightarrow \infty $. The problem is then the computation of $\lambda_n$, and since it grows with $n$, we simplify (\ref{qI}) by observing that
\begin{equation*}
\left\vert \frac{1}{e^{-2i\alpha \pi }+2h\lambda ^{\alpha }\cos (\alpha \pi)e^{-i\alpha \pi }+h^{2}\lambda ^{2\alpha }}\right\vert  \sim \frac{1}{h^{2}\lambda ^{2\alpha }}, \quad \text{for\quad }\lambda \rightarrow \infty.
\end{equation*}
Therefore, we can consider the estimate
\begin{equation} \label{gI}
q_{n}^{(I)}(\lambda )\cong \tilde{q}_{n}^{(I)}(\lambda )=\frac{4\pi \alpha e^{-\sqrt{2\alpha \overline{n}}\gamma ^{-}(\lambda )}}{h\lambda ^{\alpha }}, \text{\quad for\quad }\lambda \geq \overline{\lambda }.
\end{equation}
By imposing $\frac{d }{d\lambda}\tilde{q}_{n}^{(I)}(\lambda)=0$ and defining $s=h^{\frac{1}{\alpha}}\lambda$, after some computation, we obtain 
\begin{equation} \label{27bis}
\frac{\pi ^{2}}{((\ln s)^{2}+\pi ^{2})(\ln s+\sqrt{(\ln s)^{2}+\pi ^{2}})}=\frac{2\alpha }{\Bar{n}}.
\end{equation}%
Using $\ln s < \sqrt{(\ln s)^2+\pi^2}$, we have that
\begin{equation*}
\frac{\pi ^{2}}{((\ln s)^{2}+\pi ^{2})(\ln s+\sqrt{(\ln s)^{2}+\pi ^{2}})}\geq \frac{\pi ^{2}}{2((\ln s)^{2}+\pi ^{2})^{\frac{3}{2}}}.
\end{equation*}
This relation allows to observe that, if $s^{\star}$ is the solution of (\ref{27bis}), then there exists a constant $d$, independent of $n$, such that 
\begin{equation} \label{logaritmo}
(\ln s^{\ast })^{3}\geq dn.
\end{equation}
Now, by defining 
\begin{equation*}
\tau (s):=\frac{1}{\sqrt{1+\left(\frac{\pi }{\ln s}\right)^{2}}},
\end{equation*}
we have
\begin{equation} \label{ln(s)}
\ln s=\tau(s) \sqrt{(\ln s)^{2}+\pi ^{2}},
\end{equation}
and hence we can rewrite equation (\ref{27bis}) as
\begin{equation*}
\frac{\pi^2}{(1+\tau(s))((\ln s)^{2} + \pi^2)^{\frac{3}{2}}}= \frac{2 \alpha}{\bar{n}}.
\end{equation*}
By (\ref{logaritmo}), $\tau (s^{\ast })\sim 1$ for $n \rightarrow \infty$, and therefore
\begin{equation} \label{28bis}
s^{\ast } \sim \exp \left( \left[ \left(\frac{\Bar{n}\pi ^{2}}{4\alpha
}\right)^{\frac{2}{3}}-\pi ^{2}\right]^{\frac{1}{2}%
}\right).
\end{equation}
The above relation states that 
\begin{equation} \label{due_asterisco}
\left( \ln s^{\star} \right)^2+\pi^2 \sim \left(\frac{\bar{n} \pi^2}{4 \alpha} \right)^{\frac{2}{3}}.
\end{equation}
Let $\lambda^{\star} = h^{-1/\alpha}s^{\star}$ be the point of maximum of $\tilde{q}_n^{(I)}(\lambda)$.
Since,
\begin{equation*}
\gamma^- (\lambda) = \frac{\pi}{\sqrt{\sqrt{(\ln s)^2+\pi^2}+\ln s}},
\end{equation*}
(cf. (\ref{gamma})), using (\ref{ln(s)}) and (\ref{due_asterisco}) we have
\begin{equation*}
\gamma ^{-}(\lambda^{\ast })\sim \left( \frac{\alpha }{2\Bar{n}} \right) ^{\frac{1}{6}}\pi ^{\frac{2}{3}}.
\end{equation*}
By inserting the above relation and (\ref{28bis}) in (\ref{gI}), we finally obtain 
\begin{equation} \label{g_nI}
\max_{\lambda \geq \bar{\lambda}} q_n^{(I)}(\lambda) \cong  \tilde{q}_{n}^{(I)}(\lambda_n)\sim 4\pi \alpha e^{-c(\Bar{n}\alpha ^{2}\pi^{2})^{\frac{1}{3}}}=: g_n^{(I)},
\end{equation}
where 
\begin{equation} \label{costante_c}
c = 3 \cdot 2^{-\frac{2}{3}} \cong 1.9.
\end{equation}
\newline
$\mathbf{ q_{n}^{(II)}(\lambda )}$ -
For what concerns $q_{n}^{(II)}(\lambda )$, everything depends on the term
\begin{equation} \label{den}
\left\vert 1-e^{-i\frac{\pi }{\alpha }}h^{\frac{1}{\alpha }}\lambda
\right\vert ^{-1}.
\end{equation}
A simple analysis shows that there is a maximum at 
\begin{equation*}
h^{\frac{1}{\alpha }}\lambda =\cos \frac{\pi }{\alpha }.
\end{equation*}
Therefore for $\alpha $ such that $\cos \frac{\pi }{\alpha }\leq 0$, the function is monotone decreasing for $\lambda \geq 1$. For $\alpha $ such that $\cos \frac{\pi }{\alpha }>0$, there is a maximum for $\lambda >0$ that may be smaller or greater than $\overline{\lambda }$ (see (\ref{lb})). In any case, however, experimentally one observes that the true error is almost flat at the beginning and then follows the approximation $q_n^{(I)}(\lambda)$ (see Figure \ref{Figura1}), so that the idea is to consider the approximation 
\begin{equation} \label{gII}
\max_{1 \leq \lambda \leq \bar{\lambda}} q_{n}^{(II)}(\lambda )\cong \frac{2\pi e^{-\sqrt{2(1-\alpha
)\pi \overline{n}}}}{\sin (\alpha \pi )} =: g_{n}^{(II)},
\end{equation}%
that is obtained by neglecting the term (\ref{den}) in (\ref{qII}).
\newline

$\mathbf{ q_{n}^{(III)}(\lambda )}$ -
For $n$ large enough, the function $q_{n}^{(III)}$ is monotone decreasing (see Figure \ref{Figura2}).
Then we have 
\begin{equation} \label{3cp}
\max_{1 \leq \lambda \leq \bar{\lambda}} q_{n}^{(III)}(\lambda) = q_{n}^{(III)} (1) \cong 4 \pi \alpha h e^{-\sqrt{2(\alpha+1) \bar{n}} \gamma^+(1)},
\end{equation}
where, as before, we have neglected the term 
\begin{equation*}
\left\vert 1+2\cos (\alpha \pi )(-1)^{\alpha }h\lambda+(-1)^{2\alpha }h^{2}\lambda^2 \right\vert^{-1}
\end{equation*}
in (\ref{qIII}) to prevent the inaccuracy of our formulas whenever the poles $z_0^{(III)}$ and $z_0^{(IV)}$ belong to close parabolas.
Anyway, experimentally approximation (\ref{3cp}) is poorly accurate for small $h$ and $\alpha$.
The reason lies on the fact that the poles $\left\lbrace z_k^{(III)}\right\rbrace$ are too close to each others, and therefore formula (\ref{berr}) does not work properly. Hence, in what follows, we provide an estimate similar to (\ref{3cp}) that is obtained by removing the dependency on $h$, that is, by considering the worst case with regard to this parameter. In this view, let
\begin{equation*}
\phi (h):=\gamma^+(1) = \sqrt{\sqrt{\left( \frac{\ln h}{\alpha} \right)^2 +\pi^2} +\frac{\ln h}{\alpha}},
\end{equation*}
(see (\ref{gamma})).
It is easy to show that
\begin{equation*}
\phi(h) \sim \frac{\pi}{\sqrt{\frac{2}{\alpha}(-\ln h)}},
\end{equation*}
for $h \rightarrow 0$.
Now, by using the above approximation in (\ref{3cp}), we obtain
\begin{align}
h e^{- \sqrt{2(\alpha+1) \bar{n}}c(h)}  &\sim e^{\ln h -\sqrt{2(\alpha+1)\bar{n}}\pi \sqrt{\frac{\alpha}{2}}\frac{1}{\sqrt{-\ln h}}} \notag \\
&=e^{-y-\sqrt{p} \pi \sqrt{\frac{\alpha}{2}}\frac{1}{\sqrt{y}}}, \label{quadratino}
\end{align}
where $y=-\ln h$, $0<y<+ \infty$ for $h<1$, and $p=2(\alpha+1)\bar{n}$.
Let us consider the function 
\begin{equation*}
\xi(y)=-y-\sqrt{p}\pi \sqrt{\frac{\alpha}{2}}\frac{1}{\sqrt{y}}.
\end{equation*}
Since $\xi''(y)<0$ for $y \in (0,+\infty)$, we look for its maximum $\bar{y}$ by solving
\begin{equation*}
\frac{d}{dy}\left (-y - \sqrt{p}\sqrt{\frac{\alpha}{2}}\pi \frac{1}{\sqrt{y}}\right ) = 0 ,
\end{equation*}
that leads to 
\begin{equation*}
\bar{y} =\left( \frac{\sqrt{p}}{2} \sqrt{\frac{\alpha}{2}}\pi \right)^{\frac{2}{3}}.
\end{equation*}
Hence, we have that
\begin{equation} \label{stellina}
\xi(\bar{y})= -c p^{\frac{1}{3}} \left( \frac{\alpha}{2} \right)^{\frac{1}{3}} \pi^{%
\frac{2}{3}},
\end{equation}
where $c$ is defined in (\ref{costante_c}).
By substituting (\ref{stellina}) in (\ref{quadratino}) and going back to (\ref{3cp}), we obtain the new approximation
\begin{equation} \label{g_nIII}
\max_{1 \leq \lambda \leq \bar{\lambda}} q_n^{(III)} (\lambda) \cong 4 \pi \alpha e^{-c \left(\alpha(\alpha+1)\pi^2 \bar{n} \right)^{\frac{1}{3}}} := g_n^{(III)}.
\end{equation}
\newline

$\mathbf{ q_{n}^{(IV)}(\lambda )}$ - The behavior of the function $q_{n}^{(IV)}(\lambda )$ is very similar to the one of $q_{n}^{(II)}(\lambda )$. Therefore we consider the analog approximation
\begin{equation} \label{gIV}
\max_{\lambda \geq \bar{\lambda}} q_{n}^{(IV)}(\lambda ) \cong \frac{2\pi}{\sin(\alpha
\pi)} e^{-\sqrt{2\Bar{n}\frac{(1-\alpha)(\alpha +1)}{\alpha}\pi}}:= g_n^{(IV)}.
\end{equation}

\subsection{Comparison of bounds}

\begin{figure}
\centering
\includegraphics[scale=0.35]{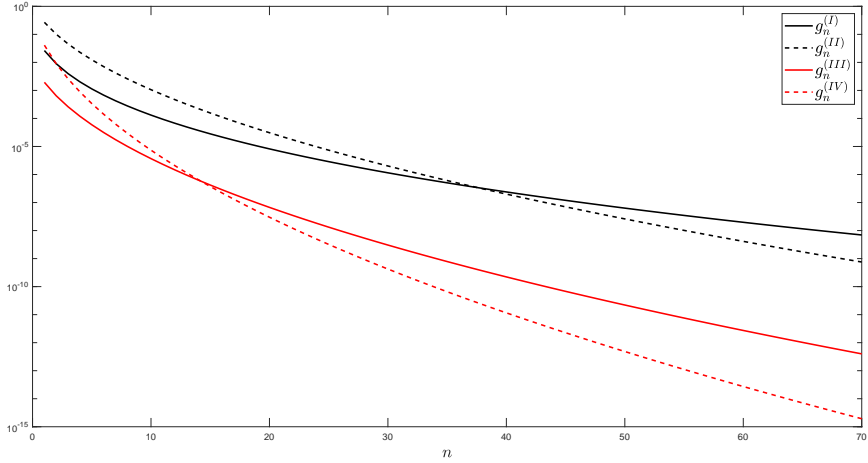}
\caption{The behavior of the functions $g_n^{(\cdot)}$ for $\alpha=0.7$ and $h=1e-2$.} \label{Figura3}
\end{figure}

By using (\ref{g_nI}), (\ref{gII}), (\ref{g_nIII}), (\ref{gIV}) we have now
\begin{align}
\Vert \mathcal{R}_{h,\alpha}(\mathcal{L} )-R_{2n-1,2n}(\mathcal{L})\Vert _{\mathcal{H}\rightarrow \mathcal{H}} & \leq \frac{\sin (\alpha \pi)}{\alpha \pi} \left(\| e_n^{(1)}(\mathcal{L}) \|_{\mathcal{H}\rightarrow \mathcal{H}}+ \| e_n^{(2)}(\mathcal{L}) \|_{\mathcal{H}\rightarrow \mathcal{H}} \right) \notag \\
&\cong \frac{\sin (\alpha \pi)}{\alpha \pi} \left(\max \left(  g_n^{(I)}, g_n^{(II)} \right) + \max \left( g_n^{(III)}, g_n^{(IV)} \right) \right). \label{bound}
\end{align} 
The next step is then the comparison of the sequences $g_n^{(\cdot)}$. We start with $g_n^{(I)}$ and $g_n^{(II)}$. Clearly $g_n^{(I)}$ is asymptotically slower and we look for $n^{\star}$ such that $g_{n}^{(I)} \geq g_{n}^{(II)}$, for $n \geq n^{\star}$.
Hence, we have to solve with respect to $n$ 
\begin{equation*}
\frac{2\pi}{\sin(\alpha\pi)} e^{-\sqrt{2\Bar{n}(1-\alpha)\pi}} = 4\pi\alpha
e^{-c(\Bar{n} \alpha^2 \pi^2)^{\frac{1}{3}}}.
\end{equation*}
By neglecting the factors before the exponentials and since $\Bar{n}= 4n+2$, we finally obtain 
\begin{equation} \label{nstar}
n^* \cong \frac{c^6}{2^5} \frac{\alpha^4}{(1-\alpha)^3} \pi - \frac{1}{2},
\end{equation}
where $c$ is defined in (\ref{costante_c}).
By (\ref{nstar}), we have $n^{\star} \geq 1$ only for $\alpha \geq \alpha^{\star} \cong 0.47$, and this is experimentally confirmed. Therefore we have
\begin{equation} \label{epsilon1}
\epsilon_n^{(1)}:= \max \left\lbrace g_n^{(I)}, g_n^{(II)} \right\rbrace =
\begin{cases}
g_n^{(I)} \quad n \geq n^{\star} \\
g_n^{(II)} \quad 1 \leq n < n^{\star}
\end{cases}.
\end{equation}
The situation is very similar for $g_n^{(III)}$ and $g_n^{(IV)}$ and we look for $n^{\star \star}$ such that $g^{(III)}_{n} \geq g^{(IV)}_{n}$, for $n \geq n^{\star\star}$, by solving
\begin{equation*}
4 \pi \alpha e^{-c \left( (\alpha+1)\bar{n} \alpha \pi^2 \right)^{\frac{1}{3}}} = \frac{2 \pi}{\sin (\alpha \pi)}e^{-\sqrt{2 \bar{n}\frac{(1-\alpha)(\alpha+1)}{2}\pi}}.
\end{equation*}
As before, by neglecting the factors before the exponentials, we find
\begin{equation} \label{nstarstar}
n^{\star \star} \cong \frac{c^6}{2^5} \frac{\alpha^5}{(1-\alpha)^3(1+\alpha)}\pi-\frac{1}{2}.
\end{equation}
Experimentally, we observe that $n^{\star \star}\geq 1$ only for $\alpha \geq \alpha^{\star \star} \cong 0.55$. Moreover, it holds $n^{\star \star}<n^{\star}$, $\forall \alpha$ (cf. (\ref{nstar}) and (\ref{nstarstar})). 
Finally, we then have
\begin{equation} \label{epsilon2}
\epsilon_n^{(2)}:=\max \left\lbrace g_n^{(III)}, g_n^{(IV)} \right\rbrace =
\begin{cases}
g_n^{(III)} \quad n \geq n^{\star \star} \\
g_n^{(IV)} \quad 1 \leq n < n^{\star \star}
\end{cases}.
\end{equation}
In Figure \ref{Figura3} we plot the sequences $g_n^{(\cdot)}$ for $\alpha=0.7$ and $h=1e-2$.
The analysis just given for $g_n^{(III)}$ and $g_n^{(IV)}$ will be important in the next section. Indeed, for what concerns the error of the method so far considered we just need to observe (see also Figure \ref{Figura3}) that 
\begin{equation*}
\max \left\lbrace g_n^{(III)}, g_n^{(IV)} \right\rbrace \ll \max \left\lbrace g_n^{(I)}, g_n^{(II)} \right\rbrace, \quad \forall n,
\end{equation*}
and in particular  that the ratio $\epsilon_n^{(2)}/\epsilon_n^{(1)}$ decays exponentially.
In this view, we finally conclude that
\begin{align} \label{stima_finale}
\Vert \mathcal{R}_{h,\alpha}(\mathcal{L} )-R_{k-1,k}(\mathcal{L})\Vert _{\mathcal{H}\rightarrow \mathcal{H}}&\cong \frac{\sin (\alpha \pi)}{\alpha \pi}  \epsilon_n^{(1)} \notag\\
&= \frac{\sin (\alpha \pi)}{\alpha \pi} 
\begin{cases}
4\pi \alpha e^{-c(\Bar{n}\alpha ^{2}\pi^{2})^{\frac{1}{3}}} \quad n \geq n^{\star} \\
\frac{2\pi e^{-\sqrt{2(1-\alpha
)\pi \overline{n}}}}{\sin (\alpha \pi )} \quad 1 \leq n < n^{\star}
\end{cases}.
\end{align}

\begin{figure} 
\begin{center}
\includegraphics[scale=0.35]{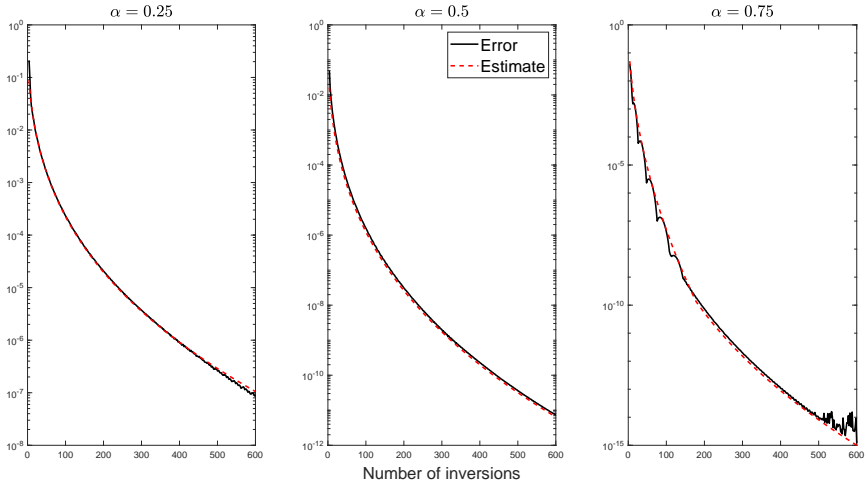}
\end{center}
\caption{Error and error estimate (\ref{stima_finale}) for the computation of $\mathcal{R}_{h,\alpha}(\mathcal{L})$, with $h=1e-2$. Different values of $\alpha$ are considered.} \label{Figura4}
\end{figure}

In order to test the behavior of the method and the accuracy of the error estimate (\ref{stima_finale}), we consider the operator
\begin{equation} \label{L}
\mathcal{L} =
{\rm diag}\left(10^0,10^{0.1},\ldots,10^{15.9}, 10^{16}\right).
\end{equation}
In Figure \ref{Figura4} we show the results for some values of $\alpha$. Here and below we always consider the spectral norm when working with matrices.
While very simple, the operator considered represents more or less the most difficult situation. Working with Matlab it is also possible to add the constant \textsf{inf} to the diagonal of $\mathcal{L}$. The results are almost indistinguishable.

\section{A balanced approach} \label{Section5}

\begin{figure} 
\begin{center}
\includegraphics[scale=0.35]{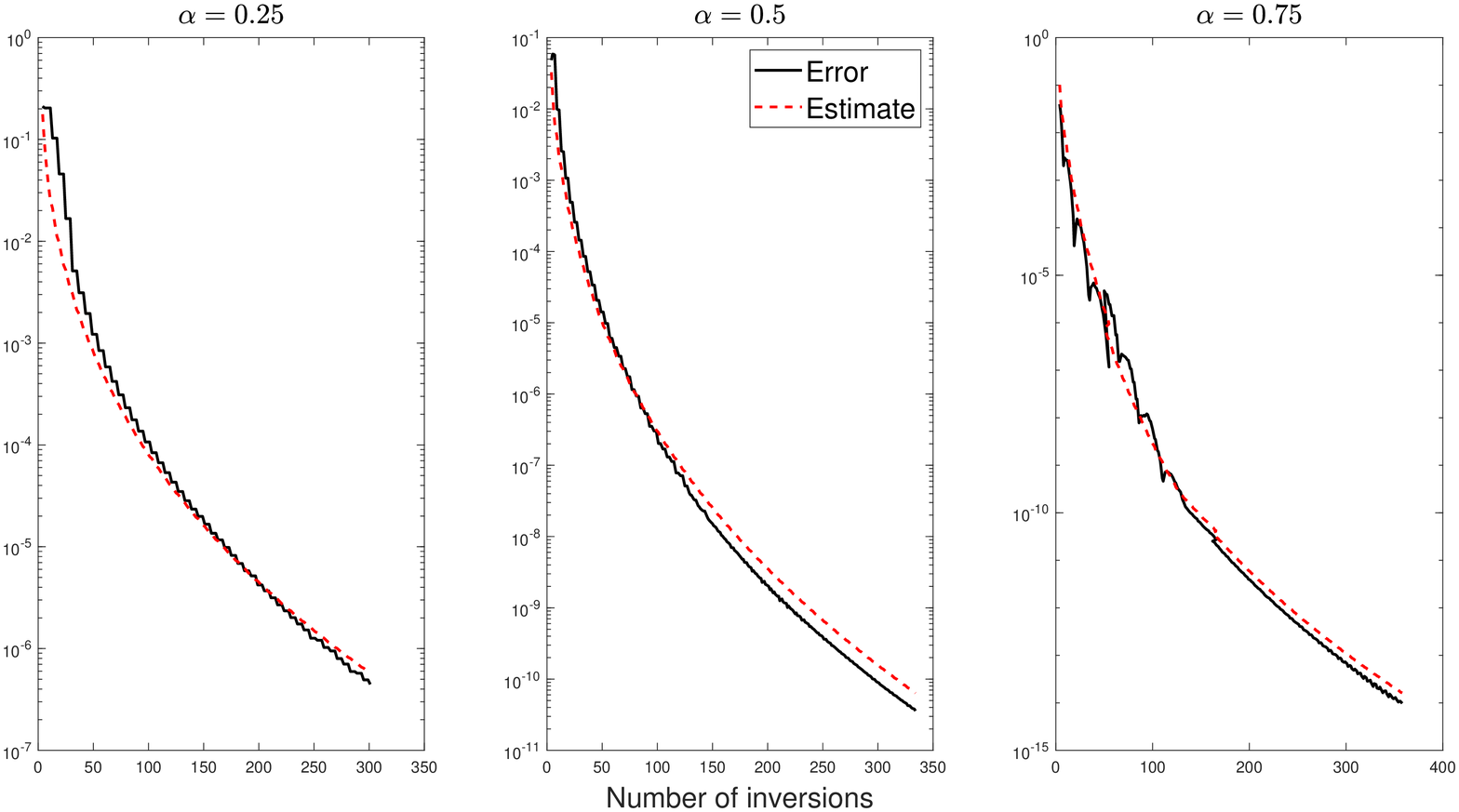}
\end{center}
\caption{The error of the balanced approach and error estimate (\ref{46bis}) for some values of $\alpha$ and $h=1e-2$.} \label{Figura5}
\end{figure}

In this section we present a modification of the Gauss-Laguerre approach that allows to reduce the number of inversions and hence the computational cost of the method, without loosing accuracy. In fact, since the computation of the first integral requires more points than the second one to achieve the same accuracy, the idea is to find $m \leq n$ such that 
\begin{equation} \label{confronto_m}
\max\left(g_n^{(I)},g_n^{(II)}\right) \cong \max \left(g_m^{(III)},g_m^{(IV)} \right),
\end{equation}
and then to consider the approximation 
\begin{equation}
\mathcal{R}_{h,\alpha}(\mathcal{L}) \cong \frac{\sin(\alpha\pi)}{\alpha\pi}
\left( I_n^{(1)} (\mathcal{L}) + I_m^{(2)} (\mathcal{L}) \right).
\end{equation} 
In this setting, we observe that the total number of inversion is $n+m$.
By using (\ref{epsilon1}), (\ref{epsilon2}) and since $n^{\star\star}\leq n^{\star}$, we have that $m$ can be obtained by imposing
\begin{equation*}
\begin{cases}
g_n^{(II)} &= g_m^{(IV)} \quad {\rm if} \; 1 \leq n \leq n^{\star\star} \\
g_n^{(II)} &= g_m^{(III)} \quad {\rm if} \; n^{\star\star} < n \leq n^{\star} \\
g_n^{(I)} &= g_m^{(III)} \quad {\rm if} \; n > n^{\star} 
\end{cases}.
\end{equation*}
After some simple computations, we find that
\begin{equation} \label{m}
m \cong
\begin{cases}
\frac{\alpha(2n+1)}{2(\alpha+1)}-\frac{1}{2} \quad {\rm if} \; 1 \leq n \leq n^{\star\star} \; {\rm and} \; n>n^{\star} \\
\frac{\left(2\sqrt{(2n+1)(1-\alpha)\pi}+\ln(2\alpha\sin(\alpha\pi))\right)^3}{27(\alpha+1)\alpha\pi^2} -\frac{1}{2} \quad {\rm if} \; n^{\star\star} < n \leq n^{\star} 
\end{cases}.
\end{equation}
An example is reported in Table \ref{Tabella1}, where $m$ is defined by using the floor operator applied to (\ref{m}). 
\begin{table} 
\caption{The values of $n$ and $m$ (\ref{m}) for $\alpha= 0.6$.} \label{Tabella1}
\begin{equation*}
\begin{array}{lccccccr}
\toprule
\mathbf{n} & 5 & 10 & 15 & 20 & 25 & 50 & 100 \\ 
\mathbf{m} & 2 & 4 & 6 & 8 & 10 & 19 & 38 \\
\bottomrule
\end{array}
\end{equation*}
\end{table}
The error is then finally estimate by (cf. (\ref{bound}))
\begin{equation} \label{46bis}
\| \mathcal{R}_{h,\alpha}(\mathcal{L})-R_{n+m-1,n+m}(\mathcal{L}) \|_{\mathcal{H}\rightarrow \mathcal{H}} \cong 2 \frac{\sin(\alpha \pi)}{\alpha 	\pi} \max \left(g_n^{(I)}, g_n^{(II)} \right).
\end{equation}
In Figure \ref{Figura5}, working with operator (\ref{L}), we plot the error and error estimate (\ref{46bis}) for some values of $\alpha$ and $h=1e-2$.
In order to have an expression of the estimates that depends on the total number of inversions $q=n+m$, by (\ref{m}) we observe that for $n>n^{\star}$
\begin{equation*}
m \cong \frac{\alpha}{\alpha+1} n,
\end{equation*}
and therefore
\begin{equation*}
q= \frac{2 \alpha +1}{\alpha+1} n.
\end{equation*}
By using (\ref{g_nI}), (\ref{costante_c}) and (\ref{46bis}) we then obtain
\begin{equation} \label{48}
\| \mathcal{R}_{h,\alpha}(\mathcal{L})-R_{q-1,q}(\mathcal{L}) \|_{\mathcal{H}\rightarrow \mathcal{H}} \cong 8 \sin(\alpha \pi) e^{-3 \left(q\frac{\alpha+1}{2\alpha+1}\alpha^2 \pi^2\right)^{\frac{1}{3}}}.
\end{equation}
Note that without balancing, that is for $q=2n$, again by (\ref{g_nI}) we have
\begin{equation*}
\| \mathcal{R}_{h,\alpha}(\mathcal{L})-R_{q-1,q}(\mathcal{L}) \|_{\mathcal{H}\rightarrow \mathcal{H}} \cong 4 \sin(\alpha \pi) e^{-3 \left(\frac{1}{2} q \alpha^2 \pi^2\right)^{\frac{1}{3}}}.
\end{equation*}
By comparing the two estimates we observe that, asymptotically, the speedup is then provided by the constant
\begin{equation*}
\frac{4}{3} < \frac{2 \alpha+2}{2\alpha+1}<2.
\end{equation*}

\section{A truncated approach} \label{Section6}

\begin{figure}
\begin{center}
\includegraphics[scale=0.35]{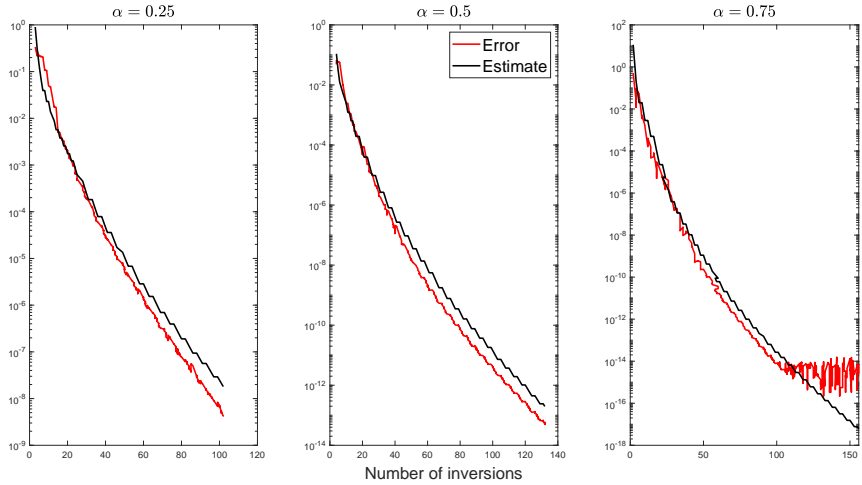}
\end{center}
\caption{The error of the truncated approach and error estimate (\ref{asterisco}) for some values of $\alpha$ and $h=1e-2$.} \label{Figura6}
\end{figure}

\begin{figure}
\begin{center}
\includegraphics[scale=0.35]{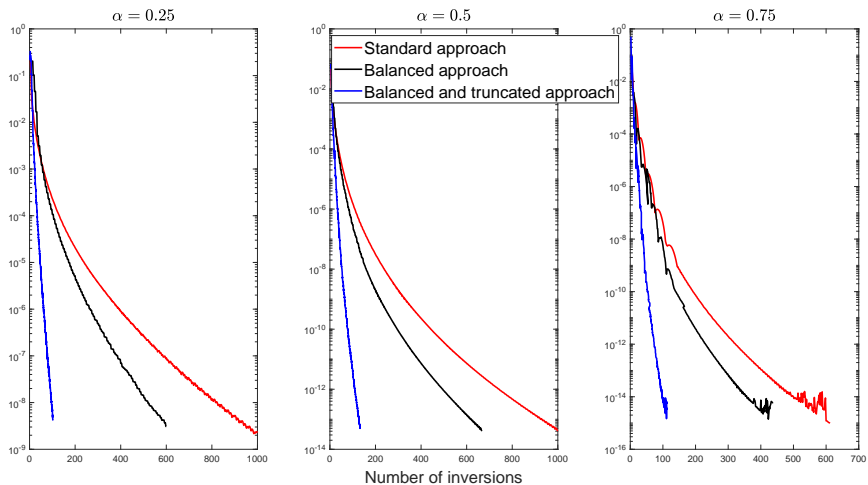}
\end{center}
\caption{Comparisons between the errors of the three approaches, standard, balanced, balanced and truncated, for different values of $\alpha$ and $h=1e-2$.} \label{Figura7}
\end{figure}

In this section we present an additional approach, already used in \cite{AN21}, to further reduce the total number of inversions without loss of accuracy. In particular, since the weights of the Gauss-Laguerre rule decay exponentially (see e.g., \cite{MO}), the idea is to find $k_n<n$ and $k_m<m$ such that we can suitably neglect the tails of the quadrature formulas $I_n^{(1)}(\mathcal{L})$ and $I_m^{(2)}(\mathcal{L})$ and therefore consider the approximation
\begin{equation} \label{approssimazione troncato}
\mathcal{R}_{h,\alpha}(\mathcal{L}) \cong \frac{\sin(\alpha \pi)}{\alpha\pi} \left(
I_{k_n}^{(1)}(\mathcal{L}) + 
I_{k_m}^{(2)}(\mathcal{L}) \right),
\end{equation}
where 
\begin{equation*}
I_{k_n}^{(1)}(\lambda)= \sum_{j=1}^{k_n} w_j^{(n)} \mathit{f}_1(x_j^{(n)}) \quad {\rm and} \quad I_{k_m}^{(2)}(\lambda) =\sum_{j=1}^{k_m}w_j^{(m)} \mathit{f}_2(x_j^{(m)}).
\end{equation*}
In this setting, the total number of inversions is $k_n+k_m$.
We start the analysis by recalling that the sequences of error approximations of the two integrals, denoted by $\{ \epsilon_n^{(1)} \}_{n\geq1}$ and $\{ \epsilon_m^{(2)} \}_{m\geq1}$ (see (\ref{epsilon1}), (\ref{epsilon2})), are such that $\epsilon_n^{(1)} \cong \epsilon_m^{(2)}$, because of the balancing introduced in previous section.
Now, remembering that, uniformly with respect to $1 \leq \lambda < + \infty$, for the functions $f_1$ and $f_2$ (see (\ref{f1}), (\ref{f2})) it holds $0 \leq f_i(x) \leq K_i$, $i=1,2$, with $K_i$ as in (\ref{K_i}), we have
\begin{equation*}
\int_0^{+\infty}e^{-x} f_i(x) dx \leq K_i \int_0^{+\infty} e^{-x}dx, \quad i=1,2.
\end{equation*}
At this point, let $s_n^{(1)}$ and $s_m^{(2)}$ be respectively the solutions of 
\begin{equation*}
K_1 \int_{s_n^{(1)}}^{+\infty} e^{-x}dx = \epsilon_n^{(1)} \quad {\rm and} \quad K_2 \int_{s_m^{(2)}}^{+\infty} e^{-x}dx = \epsilon_m^{(2)}.
\end{equation*}
From the above equations,
\begin{equation} \label{s_n}
s_n^{(1)}= -\ln\left(\frac{\epsilon_n^{(1)}}{K_1}\right) \quad {\rm and} \quad s_m^{(2)}= -\ln\left(\frac{\epsilon_m^{(2)}}{K_2}\right).
\end{equation}
Then, for the first integral we consider the truncated rule $I_{k_n}^{(1)}$, where $k_n$ is the smallest integer such that $x_j^{(n)} \geq s_n^{(1)}$, $\forall j \geq k_n$. We observe that
\begin{align*}
\left\Vert I^{(1)}(\mathcal{L}) - I_{k_n}^{(1)}(\mathcal{L})\right\Vert_{\mathcal{H}\rightarrow \mathcal{H}} & = \left\Vert I^{(1)}(\mathcal{L})
-I_n^{(1)}(\mathcal{L}) + \sum_{j=k_n + 1}^{n} w_j^{(n)}\mathit{f}_1(x_j^{(n)})
\right\Vert_{\mathcal{H}\rightarrow \mathcal{H}} \\
& \leq \epsilon_n^{(1)} + \sum_{j=k_n + 1}^{n} w_j^{(n)}\mathit{f}_1(x_j^{(n)}) \\
& \leq \epsilon_n^{(1)} + K_1\sum_{j=k_n + 1}^{n} w_j^{(n)}.
\end{align*}
Now, by using the bound (see \cite{MO})
\begin{equation*}
w_j^{(n)} \leq C(x_j^{(n)}-x_{j-1}^{(n)})e^{-x_j^{(n)}},
\end{equation*}
where C is a constant independent of $n$ and close to 1, we have (see (\ref{s_n}))
\begin{align*}
\sum_{j=k_n+1}^{n} w_j^{(n)} &\leq C \sum_{j=k_n + 1}^n
\left(x_j^{(n)}-x_{j-1}^{(n)}\right)e^{-x_j^{(n)}} \leq C \int_{x_{k_n}^{(n)}}^{+\infty} e^{-x} dx \\
&=C e^{-x_{k_n}^{(n)}} \leq C e^{-s_n^{(1)}} = C \epsilon_n^{(1)},
\end{align*}
and finally
\begin{equation} \label{1stella}
\left\Vert I^{(1)}(\mathcal{L}) -I_{k_n}^{(1)}(\mathcal{L}) \right\Vert_{\mathcal{H}\rightarrow \mathcal{H}} \leq (1+C) \epsilon_n^{(1)} \cong 2\epsilon_n^{(1)}.
\end{equation}
As for the second integral, by following the same arguments, we obtain 
\begin{equation} \label{2stelle}
\left\Vert I^{(2)}(\mathcal{L}) -I_{k_m}^{(2)}(\mathcal{L}) \right\Vert_{\mathcal{H}\rightarrow \mathcal{H}} \leq (1+C) \epsilon_m^{(2)} \cong 2 \epsilon_m^{(2)},
\end{equation}
where $k_m$ is the smallest integer such that $x_j^{(m)} \geq s_m^{(2)}$, $\forall j \geq k_m$.
It is interesting to observe that it is also possible to derive an analytical approximate expression of $k_n$ and $k_m$, that allows to understand the behavior of the error with respect to $k_n+k_m$. The analysis makes use of the relation 
\begin{equation} \label{approx_x_k}
x_{j}^{(n)} = c_{j} \frac{j^2\pi^2}{4n} \left( 1+ \mathcal{O} \left( \frac{1}{n^2} \right) \right),
\end{equation}
with $1 < c_{j} < \left( 1+\frac{1}{j} \right)^2$, given in \cite[Prop. 6.1]{AN21}.
We start with the computation of $k_n$ for $1 \leq n\leq n^*$. In this case the error $\epsilon_n^{(1)}$ is given by $g_n^{(II)}$ (see (\ref{gII}), (\ref{epsilon1})), and therefore 
\begin{equation*}
s_n^{(1)} = -\ln\left(\frac{g_n^{(II)}}{K_1}\right) = -\ln\left(\frac{2\pi}{\sin(\alpha\pi)}\right)+\sqrt{2\Bar{n}(1-\alpha)\pi},
\end{equation*}
(since $K_1=1$). 
By neglecting the term $\ln\left( \frac{2\pi}{\sin(\alpha\pi)} \right)$ and since $\Bar{n} =4n+2$, we obtain 
\begin{equation*}
s_n^{(1)} \cong \sqrt{8n(1-\alpha)\pi}.
\end{equation*}
Recalling that $k_n$ is such that $x_{k_n}^{(n)} \geq
s_n^{(1)}$ and using the relation (\ref{approx_x_k}), we try to solve with respect to $j$ 
\begin{equation*}
\sqrt{8n(1-\alpha)\pi} = c_{j} \frac{j^2\pi^2}{4n}.
\end{equation*}
By using $c_{j}\cong 1$ and the floor operator $\lfloor \cdot \rfloor$, we have that
\begin{equation} \label{jn_uno}
j_n := \left\lfloor 2(1-\alpha)^{\frac{1}{4}}\left( \frac{2n}{\pi}\right)^{\frac{3}{4}} \right\rfloor,
\end{equation}
is a good approximation of $k_n$. 
Note that $j_n \leq n$, $\forall n$, $\forall \alpha$.
From the above expression we can compute $n$ in terms of $j_n$, that is,
\begin{equation} \label{n1}
n \cong \pi\left[ \frac{j_n^4}{2^7(1-\alpha)} \right]^{\frac{1}{3}}, \quad {\rm for} \quad 1 \leq n \leq n^{\star}.
\end{equation}
Following the same steps, for $n>n^*$ we obtain 
\begin{equation} \label{jn_due}
j_n := \left\lfloor 2\sqrt{3} \left(\frac{\alpha n^2}{\pi^2}\right)^{\frac{1}{3}} \right\rfloor \cong k_n,
\end{equation}
from which 
\begin{equation} \label{n2}
n \cong \frac{\pi}{\sqrt{\alpha}} \left[ \frac{j_n}{2\sqrt{3}} \right]^{\frac{3}{2}}, \quad {\rm for} \quad n>n^*.
\end{equation}
Finally, by using the approximations (\ref{n1}) and (\ref{n2}) in (\ref{epsilon1}), we obtain
\begin{equation} \label{triangolo}
\epsilon_n^{(1)} \cong 
\begin{cases}
4 \pi \alpha e^{-3^{3/4} 2^{-1/2} \alpha^{1/2} \pi j_n^{1/2}}, \quad {\rm for} \quad n > n^{\star} \\
\frac{2 \pi}{\sin(\alpha \pi)} e^{-2^{1/3} (1-\alpha)^{1/3} \pi j_n^{2/3}}, \quad {\rm for} \quad 1 \leq n \leq n^{\star}
\end{cases}
\end{equation}
For the second integral the analysis is the same. We just need to remember that $K_2= \frac{\alpha}{\alpha+1} h^{-1/ \alpha}$.
Then, for $1 \leq m \leq m^{\star \star}$, we obtain
\begin{equation} \label{jm}
j_m := \left\lfloor \left\lbrace \frac{4m}{\pi^2} \left[ \ln \left(\frac{\alpha}{\alpha+1} h^{-\frac{1}{\alpha}} \right) + \sqrt{\frac{
8m(1-\alpha)(\alpha+1)\pi}{\alpha}}\right] \right\rbrace^{\frac{1}{2}} \right\rfloor \cong k_m,
\end{equation}
and for $m \geq m^{**}$, we have that
\begin{equation} \label{quadrato}
j_m := \left\lfloor\left\lbrace \frac{4m}{\pi^2} \left[\ln \left(\frac{\alpha}{\alpha+1} h^{-\frac{1}{\alpha}} \right) + 3\left((\alpha +1
)\alpha\pi^2m\right)^{\frac{1}{3}} \right]\right\rbrace^{\frac{1}{2}}\right\rfloor \cong k_m.
\end{equation}
At this point we are able to write down the final error estimates. From (\ref{1stella}), (\ref{2stelle}), (\ref{triangolo}) and since $\epsilon_n^{(1)} \cong \epsilon_m^{(2)}$, we obtain
\begin{align}
&\left\Vert \mathcal{R}_{h,\alpha} (\mathcal{L})- R_{q-1,q}(\mathcal{L}) \right\Vert_{\mathcal{H} \rightarrow \mathcal{H}} \notag\\
&\leq \frac{\sin (\alpha \pi)}{\alpha \pi} \left( \left\Vert I^{(1)}(\mathcal{L})-I_{j_n}^{(1)}(\mathcal{L}) \right\Vert_{\mathcal{H} \rightarrow \mathcal{H}} + \left\Vert I^{(2)}(\mathcal{L})-I_{j_m}^{(2)}(\mathcal{L}) \right\Vert_{\mathcal{H} \rightarrow \mathcal{H}} \right) \notag \\
&\cong 4 \frac{\sin(\alpha \pi)}{\alpha \pi} \epsilon_n^{(1)} \notag \\
&\cong 4 \frac{ \sin(\alpha \pi)}{\alpha \pi}
\begin{cases}
4 \pi \alpha e^{-c \pi 2^{1/6}3^{-1/4}\alpha^{1/2} j_n^{1/2}}, \quad {\rm for} \quad n > n^{\star} \\
\frac{2 \pi}{\sin(\alpha \pi)} e^{-3^{3/4} 2^{-1/2} \alpha^{1/2} \pi j_n^{1/2}}, \quad {\rm for} \quad 1 \leq n \leq n^{\star}
\end{cases}, \label{asterisco}
\end{align} 
where $q=j_n+j_m$.
In Table \ref{Tabella2} we show the values of $m$, $k_n$, $k_m$, together with the theoretical approximations $j_n$ and $j_m$, with respect to $n$, for the case of $\alpha=0.75$. It is rather clear that the approximations provided by $j_n$ and $j_m$ are fairly accurate.
\begin{table} 
\caption{The values of $m$, $k_n$, $j_n$, $k_m$ and $j_m$
with respect to $n$ in the case of $\alpha= 0.75$.} \label{Tabella2}
\begin{equation*}
\begin{array}{cccccccc}
\toprule
\mathbf{n } & 5 & 10 & 15 & 20 & 25 & 50 & 100 \\ 
\mathbf{m } & 2 & 4 & 7 & 9 & 11 & 16 & 46 \\ 
\mathbf{k_n  (j_n)} & 3(2) & 5(4) & 7(6) & 9(8) & 10(10) & 18(18) & 30(30) \\ 
\mathbf{k_m  (j_m)} & 2(2) & 4(4) & 5(6) & 6(6) & 7(8) & 11(10) & 21(22) \\ 
\bottomrule
\end{array}
\end{equation*}
\end{table}
In order to have an asymptotic expression of estimate (\ref{asterisco}) that depends on the total number of inversions $q$, from (\ref{quadrato}) we first consider the approximation
\begin{equation*}
j_m \cong [12 \pi^{-\frac{4}{3}} (\alpha+1)^{\frac{1}{3}} \alpha^{\frac{1}{3}} m^{\frac{4}{3}} ]^{\frac{1}{2}}, \quad {\rm for} \quad m > m^{\star},
\end{equation*}
and hence
\begin{equation*}
m \cong \frac{\pi}{2^{3/2} 3^{3/4} (\alpha+1)^{1/4} \alpha^{1/4}} j_m^{\frac{3}{2}}.
\end{equation*}
By using the above approximation in $g_m^{(III)}$ (see (\ref{g_nIII})), we obtain
\begin{equation*}
g_m^{(III)} \cong 4 \pi \alpha e^{- 3^{3/4} 2^{1/2} \alpha^{1/4} (\alpha+1)^{1/4} \pi j_m^{1/2}}.
\end{equation*}
Moreover, since $\epsilon_n^{(1)} \cong \epsilon_m^{(2)}$, we have that 
\begin{equation*}
3^{\frac{3}{4}} 2^{-\frac{1}{2}} \pi \alpha^{\frac{1}{2}} j_n^{\frac{1}{2}} \cong 3^{\frac{3}{4}} 2^{-\frac{1}{2}} \pi \alpha^{\frac{1}{4}} (\alpha+1)^{\frac{1}{4}} j_m^{\frac{1}{2}},
\end{equation*}
and therefore
\begin{equation*}
j_n+j_m \cong \left[ 1 +\left( \frac{\alpha}{\alpha+1} \right)^{\frac{1}{2}}\right] j_n.
\end{equation*}
Then, by (\ref{asterisco}), we finally have that asymptotically
\begin{equation} \label{asintotica} 
\left\Vert \mathcal{R}_{h,\alpha} (\mathcal{L})- R_{q-1,q}(\mathcal{L}) \right\Vert_{\mathcal{H} \rightarrow \mathcal{H}} \cong 16 \sin (\alpha \pi) \exp \left(-3^{\frac{3}{4}} 2^{-\frac{1}{2}} \pi \alpha^{\frac{1}{2}} \left[1+ \left( \frac{\alpha}{\alpha+1} \right)^{\frac{1}{2}} \right]^{-\frac{1}{2}} q^{\frac{1}{2}} \right).
\end{equation}
As an example of the remarkable improvements of the balanced and truncated approach, working with operator (\ref{L}), in Figure \ref{Figura6} we plot the error and error estimate (\ref{asterisco}), while in Figure \ref{Figura7} we compare the three approaches developed in this work, for different values of $\alpha$ and $h=1e-2$.
Finally, in Algorithm \ref{Algoritmo} we summarize the steps necessary to implement the method.
\begin{algorithm}[Balanced and truncated Laguerre method] \label{Algoritmo}
Input $\alpha, h ,\mathcal{L}$

\noindent evaluate $n^{\star}$, $n^{\star \star}$ using (\ref{nstar}), (\ref{nstarstar})

\medskip
\noindent for $n=1,\ldots$
\begin{itemize}
\item [] evaluate $m$ using (\ref{m}) 
\item [] compute $w_j^{(n)}$, $x_j^{(n)}$ and $w_j^{(m)}$, $x_j^{(m)}$ 
\item [] evaluate $j_n$ ((\ref{jn_uno}), (\ref{jn_due})) and $j_m$ ((\ref{jm}), (\ref{quadrato})) 
\item [] calculate approximation (\ref{approssimazione troncato})
\end{itemize}
\end{algorithm}

\section{Conclusions}

In this work we have described an efficient method for the computation of the resolvent of the fractional powers in the continuous setting of a generic Hilbert space. 
The use of the Gauss-Laguerre rule, with the improvements developed in Section \ref{Section5} and \ref{Section6}, leads to a method whose rate of convergence is the same of the scalar case, that is of type $\exp(-\emph{const} \cdot q^{1/2})$, where $q$ represents the number of inversions (cf. the definitions of $q_n^{(\cdot)}(\lambda)$ in Section \ref{Section3.2} and formula (\ref{asintotica})).
Moreover, we have provided accurate error estimates even if in Section \ref{Section4.1} we have been forced to adopt approximations only justified by experimental evidences.
We also remark that the final algorithm (Algorithm \ref{Algoritmo}) does not require the definition of any parameter. It only needs the code for the computation of the Laguerre nodes and weights, for which we have employed the Matlab function \textsf{lagpts.m} from \textsf{chebfun} (see \cite{Tbook}).

\section*{Acknowledgements}

This work was partially supported by GNCS-INdAM, FRA-University of Trieste and CINECA under HPC-TRES program award number 2019-04. Eleonora Denich and Paolo Novati are members of the INdAM research group GNCS.

\end{document}